\documentclass[12pt,dvips]{amsart}
\usepackage[usenames]{color}
\usepackage{tikz}

\newtheorem{lemma}{Lemma}[section]
\newtheorem{theorem}{Theorem}[section]
\newtheorem{definition}{Definition}[section]
\newtheorem{proposition}{Proposition}[section]

\DeclareMathOperator{\sign}{sign}

\DeclareMathOperator{\CPP}{CPP}
\DeclareMathOperator{\RPP}{RPP}

\DeclareMathOperator{\lhs}{lhs}

\DeclareMathOperator{\PP}{PP}

\title{Enumeration of Cylindric Plane Partitions - part I}
\author{Robin Langer}
\date{\today}

\begin{document}
\newdimen{\cellsize}
\newcommand\Bigboxes{\setlength{\cellsize}{24pt}\def\boxformat{}}%
\newcommand\bigboxes{\setlength{\cellsize}{18pt}\def\boxformat{}}
\newcommand\medboxes{\setlength{\cellsize}{14pt}\def\boxformat{}}
\newcommand\smallboxes{\setlength{\cellsize}{8pt}\def\boxformat{\scriptstyle}}
\medboxes
\newsavebox{\cellcontent}
\def\hidehrule#1#2{\kern-#1
  \hrule height#1 depth#2 \kern-#2 }%
\def\hidevrule#1#2{\kern-#1{\dimen\cellcontent=#1%
    \advance\dimen\cellcontent by#2\vrule width\dimen\cellcontent}\kern-#2 }%
\def\makeblankbox#1#2{\hbox{\lower\dp\cellcontent\vbox{\hidehrule{#1}{#2}%
    \kern-#1 
    \hbox to \wd\cellcontent{\hidevrule{#1}{#2}%
      \raise\ht\cellcontent\vbox to #1{}
      \lower\dp\cellcontent\vtop to #1{}
      \hfil\hidevrule{#2}{#1}}%
    \kern-#1\hidehrule{#2}{#1}}}}
\newcommand\cellify[1]{\defaultcell%
\sbox{\cellcontent}{\vbox to \cellsize{%
\vfill%
\hbox to \cellsize{\hfill$\boxformat #1$\hfill}
\vfill}}%
\rlap{\drawnbox}
\usebox{\cellcontent}}
\newcommand\tableau[1]{\vtop{\let\\\cr
\baselineskip -16000pt \lineskiplimit 16000pt \lineskip 0pt
\ialign{&\cellify{##}\cr#1\crcr}}}
\newcommand\defaultcell{\gdef\drawnbox{
\makeblankbox{0.2pt}{0.2pt}
}}
\newcommand\graycell{\gdef\drawnbox{%
\rlap{\color{Gray}\vrule width \cellsize height \cellsize}%
\makeblankbox{0.2pt}{0.2pt}
}}
\newcommand\bluecell{\gdef\drawnbox{%
\rlap{\color{Blue}\vrule width \cellsize height \cellsize}%
\makeblankbox{0.2pt}{0.2pt}
}}
\newcommand\redcell{\gdef\drawnbox{%
\rlap{\color{Red}\vrule width \cellsize height \cellsize}%
\makeblankbox{0.2pt}{0.2pt}
}}
\newcommand\greencell{\gdef\drawnbox{%
\rlap{\color{Green}\vrule width \cellsize height \cellsize}%
\makeblankbox{0.2pt}{0.2pt}
}}
\newcommand\orangecell{\gdef\drawnbox{%
\rlap{\color{Orange}\vrule width \cellsize height \cellsize}%
\makeblankbox{0.2pt}{0.2pt}
}}
\newcommand\yellowcell{\gdef\drawnbox{%
\rlap{\color{Yellow}\vrule width \cellsize height \cellsize}%
\makeblankbox{0.2pt}{0.2pt}
}}
\newcommand\thickcell{\gdef\drawnbox{
\makeblankbox{0.2pt}{0.1\cellsize}%
}}
\newcommand\missingcell{\gdef\drawnbox{}}
\newcommand\vdotscell{\gdef\drawnbox{\kern-1.6pt\vbox{\baselineskip=4pt\lineskiplimit=0pt\hbox{}\hbox{.}\hbox{.}\hbox{.}\hbox{}}}}
\newcommand\hdotscell{\gdef\drawnbox{\vbox to \cellsize{\hbox{\kern1pt$\ldotp\ldotp\ldotp$}}}}
\newcommand\vhdotscell{\gdef\drawnbox{\rlap{\kern-1.6pt\vbox{\baselineskip=4pt\lineskiplimit=0pt\hbox{}\hbox{.}\hbox{.}\hbox{.}\hbox{}}}\vbox to \cellsize{\hbox{\kern1pt$\ldotp\ldotp\ldotp$}}}}
%
\newcommand\vertlinecell{\gdef\drawnbox{\unitlength=\cellsize%
\begin{picture}(1,1)
\put(0,0){\line(0,1){1}}
\end{picture}}}
\newcommand\horizlinecell{\gdef\drawnbox{\unitlength=\cellsize%
\begin{picture}(1,1)
\put(0,1){\line(1,0){1}}
\end{picture}}}
\newcommand\defaultcella{\gdef\drawnbox{
\unitlength=\cellsize%
\begin{picture}(1,1)
\put(0,0){\line(1,0){1}}
\put(0,0){\line(0,1){1}}
\put(1,0){\line(0,1){1}}
\put(0,1){\line(1,0){1}}
\end{picture}}}
\newcommand\thickcella{
\gdef\drawnbox{%
\unitlength=\cellsize%
\begin{picture}(1,1)
\linethickness{0.1\cellsize}
\put(0.0,0.05){\line(1,0){1}}
\put(0.05,0){\line(0,1){1}}
\put(0.95,0){\line(0,1){1}}
\put(0,0.96){\line(1,0){1}}
\end{picture}}}
\newcommand\graycella{\gdef\drawnbox{
\rlap{\color{Gray}\vrule width \cellsize height \cellsize}%
\unitlength=\cellsize%
\begin{picture}(1,1)
\put(0,0){\line(1,0){1}}
\put(0,0){\line(0,1){1}}
\put(1,0){\line(0,1){1}}
\put(0,1){\line(1,0){1}}
\end{picture}}%
}

%
\mbox{}
\vspace{2ex}

\maketitle

\begin{abstract}
Cylindric plane partitions may be thought of as a natural generalization of reverse plane partitions. A generating series for the enumeration of cylindric plane partitions was recently given by Borodin. The first result of this paper is a $(q,t)$-analog of Borodin's identity which extends previous work by Okada in the reverse plane partition case. Our proof uses commutation relations for $(q,t)$-vertex operators acting on Macdonald polynomials as given by Garsia, Haiman and Tesla. The second result of this paper is an explicit combinatorial interpreation of the $(q,t)$-Macdonald weight in terms of a non-intersecting lattice path model on the cylinder.
\end{abstract}

\section{Introduction}

Cylindric plane partitions were first introduced by Gessel and Krattenthaler \cite{gessel-1997}. 
We shall work with a modified, though equivalent, definition.

\begin{definition} \label{cylindric-def}
For any binary string $\pi$ of length $T$,
a \emph{cylindric plane partition} with profile $\pi$ may be defined as a sequence of integer partitions:
\begin{equation} 
(\mu^0, \mu^1, \ldots \mu^T) \qquad \qquad \mu^0 = \mu^T 
\end{equation}
such that if $\pi_k = 1$ then $\mu^k / \mu^{k-1}$ is a \emph{horizontal strip}, otherwise if $\pi_k = 0$ then $\mu^{k-1} / \mu^k$ is a horizontal strip.
\end{definition}

For the precise definition of an integer partition and a horizontal strip, see section \ref{definitions}.
In the special case where $\mu^0 = \mu^T = \emptyset$ we recover the usual definition of a reverse plane partition (see, for example \cite{adachi} for a nice review).
If, in addition to this there are no inversions in the profile (see definition \ref{inversion}) then we have a regular plane partition.

A ``cube'' of a cylindric plane partition is defined to be a ``box'' of one of the underlying integer partitions.

\begin{definition}
The \emph{weight} of the cylindric partition $\mathfrak{c} = (\mu^0, \mu^1, \ldots \mu^T)$ is given by
$  |\mathfrak{c}| = |\mu^1| + |\mu^2| + \cdots |\mu^T|$.
\end{definition}

In other words, the weight of a cylindric plane partition is the number of cubes. Note that to avoid double counting, we do not include the boxes of the partition $\mu^0$ in the definition of the weight of $\mathfrak{c}$.

The beginning of the enumerative theory of plane partitions is the following famous identity of MacMahon \cite{macmahon}:
\begin{equation} \label{macmahon}
\sum_{\mathfrak{c} \in \PP} z^{|\mathfrak{c}|} = \left ( \frac{1}{1-z^n} \right )^n 
\end{equation}

The sum on the left hand side is over all regular plane partitions. If, in addition to this there are no inversions in the profile (see section \ref{definitions}) then we have a regular plane partition.

MacMahon's original proof involved delicate combinatorial arguments. It was Okounkov \cite{okounkov} who first pointed out that enumerative results for plane partitions may be obtained by considering commutation relations between vertex operators acting on fermionic fock space \cite{soliton}. The underlying algebraic structure is that of the Heisenberg algebra. By the boson-fermion correspondence these operators may be alternatively thought of as acting on symmetric functions.
The \emph{Pieri rules} for Schur functions are key to this approach:
\begin{equation}\label{p1} S_\mu[X] h_r[X] = \sum_{\lambda \in U_r(\mu)} S_\lambda[X] \end{equation}
\begin{equation}\label{p2} S_\lambda[X+z] = \sum_r \sum_{\mu \in D_r(\mu)} S_\mu[X] z^r \end{equation}

Here $U_r(\mu)$ denotes the set of all partitions which can be obtained from $\mu$ by adding a horrizontal $r$-strip and $D_r(\lambda$) denotes the set of all partitions which can be obtained from $\lambda$ by removing a horizontal $r$-strip.

The next important result in the subject is the following hook-product formula for the enumeration of reverse plane partitions with arbitrary profile $\pi$ which is due to Stanley \cite{stanley}:
\begin{equation} \label{stanley}
\sum_{\mathfrak{c} \in \RPP(\pi)} z^{|\mathfrak{c}|} = 
\prod_{\substack{i < j \\ \pi_i > \pi_j}} \frac{1}{1-z^{j-i}} 
\end{equation}

More recently the following hook-product formula for the enumeration of cylindric plane partitions of given profile $\pi$ of length $T$ was given by Borodin \cite{borodin}: 
\begin{equation} \label{borodin}
\sum_{\mathfrak{c} \in \CPP(\pi)} z^{|\mathfrak{c}|} = 
\prod_{n \geq 0} \left ( \frac{1}{1-z^{nT}} \prod_{\substack{i < j \\ \pi_i > \pi_j}} \frac{1}{1-z^{j-i + nT}} 
\prod_{\substack{i > j \\ \pi_i > \pi_j }} \frac{1}{1 - z^{j-i + (n+1)T}} \right ) 
\end{equation}

Borodin's proof uses the same vertex operator idea as Okounkov. A very different proof involving the representation theory of $\widehat{sl_n}$ was later given by Tingley \cite{tingley}

Macdonald polynomials are a natural $(q,t)$-deformation of the classical Schur polynomials. The \emph{Pieri rules} for Macdonald polynomials are very similar to those for the Schur polynomials, only in the Macdonald case certain coefficients appear (\cite{macdonald} page 341).

\begin{equation}\label{pieri-coeff}
 \varphi_{\lambda / \mu}(q,t) = \prod_{s \in C_{\lambda / \mu}}\frac{1 - q^{a_\lambda(s) + 1} t^{\ell_\lambda(s)}}{1 -q^{a_\lambda(s)} t^{\ell_\lambda(s)+1}} 
\prod_{s \in {C}_{\lambda / \mu}}\frac{1 - q^{a_\mu(s)} t^{\ell_\mu(s)+1}}{1 -q^{a_\mu(s)+1} t^{\ell_\mu(s)}}
\end{equation}
\begin{equation}\label{pieri-coeff2}
 \psi_{\lambda / \mu}(q,t) = \prod_{s \not\in C_{\lambda / \mu}}\frac{1 - q^{a_\lambda(s)+1} t^{\ell_\lambda(s)}}{1 - q^{a_\lambda(s)} t^{\ell_\lambda(s)+1}} 
\prod_{s \not\in {C}_{\lambda / \mu}}\frac{1 -q^{a_\mu(s)} t^{\ell_\mu(s)+1}}{1 -q^{a_\mu(s)+1} t^{\ell_\mu(s)}}
\end{equation}

Here $C_{\lambda / \mu}$ denotes the set of columns of $\lambda$ which are longer than the corresponding columns of $\mu$. For any box $s$ define $a_\lambda(s)$ to be the ``arm length'' of $s$ and $\ell_\lambda(s)$ to be the ``leg length'' of $s$ (see section \ref{definitions}).

By using Macdonald polynomials instead of Schur functions, Okada \cite{okada} obtained the following $(q,t)$-deformation of Stanley's result:

\begin{equation} \sum_{\mathfrak{c} \in \RPP(\pi)} W_\mathfrak{c}(q,t) z^{|\mathfrak{c}|} = 
\prod_{\substack{i < j \\ \pi_i > \pi_j}} \frac{(tz^{j - i};q)_\infty}{(z^{j - i};q)_\infty}  
\end{equation}
Note that we are making use of the hypergoemetric notation:
\begin{equation}
(a;q)_\infty = \prod_{n \geq 0} (1 - aq^n)
\end{equation}

If $\mathfrak{c} = (\mu^0,\mu^1, \ldots \mu^T)$ then the weight function is given by:

\begin{equation}  \label{weight}
W_{\mathfrak{c}}(q,t) = \prod_{\substack{k =1 \\ \pi_k = 1}}^T \varphi_{\mu^k / \mu^{k-1}}(q,t) \prod_{\substack{k = 1\\ \pi_k = 0}}^T \psi_{\mu^{k-1} / \mu^k}(q,t) 
\end{equation}

Observe that when $q=t$, Okada's formula reduces to that of Stanley. The regular plane partition case of Okada's identity had been previously given by Vuletic \cite{vuletic2}. 

The first result of this paper is an analogous $(q,t)$-deformation of Borodin and Tingley's formula for the enumeration of cylindric plane partitions:

\begin{theorem}\label{qt-borodin}
\begin{equation} \label{maintheorem}
\sum_{\mathfrak{c} \in CPP(\pi)} W_\mathfrak{c}(q,t) z^{|\mathfrak{c}|} = 
\prod_{n \geq 0} \left ( \frac{1}{1 - z^{nT}}  
\prod_{\substack{i < j \\ \pi_i > \pi_j}} \frac{(tz^{j - i + nT};q)_\infty}{(z^{j - i + nT};q)_\infty} 
\prod_{\substack{i > j \\ \pi_i > \pi_j}} \frac{(tz^{j-i + (n+1)T};q)_\infty}{(z^{j-i + (n+1)T};q)_\infty} \right )
\end{equation}
\end{theorem}

The weight function is exactly the same as that given by Okada (equation \ref{weight}).
When $q=t$ one finds that Theorem (\ref{qt-borodin}) reduces to equation (\ref{borodin}). 
The Hall--Littlewood case ($q=0$) of Theorem (\ref{qt-borodin}) has been previously given in Corteel, Savelief and Vuletic \cite{vuletic}.

The proof of Theorem  (\ref{qt-borodin}) uses commutation relations for certain $(q,t)$-vertex operators acting on Macdonald polynomials which are essentially due to Garsia, Haiman and Tesla \cite{garsia}. 

The nature of the proof is such that the identity remains true if on the left hand side we replace:
\[z^{\mathfrak{|c|}} \mapsto z_0^{|\mu_0|} z_1^{|\mu_1|} \cdots z_{T-1}^{|\mu_{T-1}|}\]
while on the right hand side we replace:
\begin{align*} 
z^{nT} & \mapsto z_0^{n} z_1^{n} \cdots z_{T-1}^{n} \\
z^{j - i + nT} & \mapsto z_0^{n} z_1^n \cdots z_i^{n} z_{i+1}^{n+1} \cdots z_{j}^{n+1} z_{j+1}^n \cdots z_{T-1}^n \quad \quad \quad \quad \text{ when } i < j \\
z^{j-i +(n+1)T} & \mapsto z_0^{n+1} z_1^{n+1} \cdots z_j^{n+1} z_{j+1}^n + \cdots z_i^{n} z_{i+1}^{n+1} \cdots z_{T-1}^{n+1} \quad \text{ when } i > j
\end{align*}

This provides a new refined version of Borodin's identity even in the Schur case, though the refined version of the reverse plane partition case had been previously given by Okada.

Recall that in the plethystic notation \cite{garsia}, if $a(q,t) = \sum_{n,m} a_{n,m} \, q^n t^m $
with $a_{n,m} \in \mathbb{Z} $ and $a_{0,0} = 0$, then we have:
\begin{equation} \Omega \left [ a(q,t) \right ] = \prod_{n,m} \frac{1 }{(1 - q^n t^m)^{a_{n,m}}} \end{equation}

Making use of this notation, the cylindric weight function may be given an explicit combinatorial description:

\begin{theorem}\label{comb-weight}
\begin{equation}\label{cylindric-weight}
W_\mathfrak{c}(q,t) = \Omega \left [ (q-t) \mathcal{D}_\mathfrak{c}(q,t) \right ] 
\end{equation}
where the alphabet $\mathcal{D}_\mathfrak{c}(q,t)$ is given by:
\begin{equation}\mathcal{D}_\mathfrak{c}(q,t) =  \sum_{s \in peak(\mathfrak{c})} \, q^{a_\mathfrak{c}(s)} t^{\ell_\mathfrak{c}(s)} - \sum_{s \in valley(\mathfrak{c})} \, q^{a_\mathfrak{c}(s)} t^{\ell_\mathfrak{c}(s)}
\end{equation}
\end{theorem}

The precise definition of ``valley'' and ``peak'' cubes will be given in section \ref{lattice-paths}. Theorem \ref{comb-weight} reduces to the combinatorial formula for the Hall--Littlewood weight function in the plane partition case in \cite{vuletic2} and both the reverse plane partition and the cylindric plane partition case in \cite{vuletic}.

The outline of this paper is as follows. In section \ref{definitions} we clarify a number of definitions pertaining to integer partitions.
In section \ref{lattice-paths} we introduce a model of non-intersecting lattice paths on a cylinder.
In section \ref{symfun} we recall some basic results from the theory of symmetric functions and Macdonald polynomials.
In section \ref{theorem1} we prove Theorem \ref{qt-borodin}.
In section \ref{theorem2} we prove Theorem \ref{comb-weight}.
Finally in section \ref{conclusion} we suggest some possible avenues for future research.

\section{Definitions}\label{definitions}

An \emph{integer partition} is simply a weakly decreasing list of non-negative integers which eventually stabilizes at zero. 
If the sum of the parts of $\lambda$ is equal to $n$, then we say that $\lambda$ is a partition of $n$ and write $|\lambda| = n$. 
The conjugate of the integer partition $\lambda = (\lambda_1, \lambda_2, \ldots, \lambda_k)$ is defined to be $\lambda' = (\lambda'_1, \lambda'_2, \ldots \lambda'_r)$ where $\lambda'_j = \#\{i \,|\, \lambda_i \geq j \}$.

It is often convenient to represent an integer partition visually as a \emph{Young diagram}, which is a collection of boxes in the cartesian plane which are ``stacked up'' in the bottom right hand corner. Note that our convention differs from both the standard French and English conventions. 
The \emph{minimum profile} of an integer partition is the binary string which traces out the ``jagged boundary'' of the associated young diagram. Reading from the top right hand corner to the bottom left hand corner, a zero is recorded for every vertical step and a one for every horizontal step. For example the minimum profile of our example partition $\lambda = (5,3,3,2)$ is $110100110$:
\[
\tableau{
\missingcell & \missingcell & \missingcell & \missingcell 1 & \missingcell 1 \\
\missingcell & \missingcell & \missingcell 1 & 0 & \\
\missingcell & \missingcell & 0 & &  \\
\missingcell 1 & \missingcell1 & 0 & &  \\
0 & &&&\\
} 
\]

The minimum profile of an integer partition necessarily starts with a one and ends with a zero.  An integer partition is uniquely determined by its minimum profile.
A \emph{generalized profile} is an arbitrary string of zeros and ones.
Each generalized profile associated to a minimum profile, and hence an integer partition, by removing the leading zeros and trailing ones.

\begin{definition}\label{inversion}
An \emph{inversion} in a binary string $\pi$ is a pair of indices $(i,j)$ such that $i < j$ and $\pi_i > \pi_j$. 
\end{definition}
There is a natural bijection between the ``boxes'' of an integer partition $\lambda$ and the inversions in any generalized profile of $\lambda$.
In order to work with the model of non-intersecting lattice paths on the cylinder it is necessary to pass from ``cartesian coordinates'' to ``inversion coordinates''.

The box $s \in \lambda$ with ``cartesian coordinates'' $(i,j)$ 
has \emph{arm length}  given by $a_\lambda(s) = \lambda_i - j$
and \emph{leg length} given by $\ell_\lambda(s) = \lambda'_j - i$. 
The hook length of the box $s$ is defined to be $h_\lambda(s) = a_\lambda(s) + b_\lambda(s) + 1$.

If the box $s$ has ``inversion coordinates'' $(i,j)$ then the arm length is given by $a_\lambda(s) = \#\{ i < k < j \, | \, \pi_k = 1\}$
and the leg length is given by $\ell_\lambda(s) = \#\{ i < k < j \, | \, \pi_k = 0\}$

We say that $\mu \subseteq \lambda$ if and only if $\mu_i \leq \lambda_i$ for all $i$.
For any pair of partitions $\lambda$ and $\mu$ satisfying $\mu \subseteq \lambda$ we say that $\lambda / \mu$ is a \emph{horizontal strip} and write $\mu \preceq \lambda$ if the following interlacing condition is satisfied:
\[ \lambda_1 \geq \mu_1 \geq \lambda_2 \geq \mu_2 \cdots \]


The following two lemmas are a straightforeward consequence of the definitions, they will nevertheless be needed
for the bijection between cylindric plane partitions and non-intersecting lattice paths on the cylinder.

\begin{lemma}\label{hstripcol}
$\lambda / \mu $ is a horizontal strip if and only if for each $j$ we have: 
\[ \lambda'_j - \mu'_j \in \{0,1\} \]
\end{lemma}

\begin{lemma}\label{hstripbin}
$\lambda / \mu$ is a horizontal strip if and only if the difference between the position of the $k$-th one in the profile of $\lambda$ and the position of the $k$-th one in the profile of $\mu$ is equal to zero or one.
\end{lemma}


\section{Lattice Paths}\label{lattice-paths}

The goal of this section is to give a bijection between cylindric plane partitions, defined as periodic interlacing sequences, and certain families of non-intersecting lattice paths on the cylinder - 
or equivalently rhombus tilings on the cylinder. Although this bijection is well known for reverse plane partitions,
we could not find it stated explicitly in the literature for the cylindric case.

Before proceeding any further, here is an example:

\begin{center}
\begin{tikzpicture}[scale=0.8]

\begin{scope}
\draw[step=0.5cm,color=gray] (0,0) grid (2.5,11);

\foreach \i in {0,...,2}
{
	\foreach \j in {0,...,11}
	{
		\path (\i,\j) node[fill,shape=circle,inner sep=0.3mm](A\i\j){};
	}
}

\foreach \i in {0,...,2}
{
	\foreach \j in {0,...,10}
	{
		\path (\i+0.5,\j+0.5) node[fill,shape=circle,inner sep=0.3mm](B\i\j){};
	}
}

\draw[thick] (A01) -- (B01) -- (A11) -- (B11) -- (A21) -- (B21);
\draw[thick] (A02) -- (B02) -- (A12) -- (B12) -- (A23) -- (B22);
\draw[thick] (A05) -- (B04) -- (A14) -- (B14) -- (A25) -- (B25);
\draw[thick] (A07) -- (B06) -- (A16) -- (B16) -- (A27) -- (B27);
\draw[thick] (A08) -- (B07) -- (A18) -- (B18) -- (A29) -- (B28);
\draw[thick] (A09) -- (B09) -- (A19) -- (B19) -- (A210) -- (B29);
\draw[thick] (A010) -- (B010) -- (A111) -- (B110) -- (A211) -- (B210);

\path (0,0) node[draw,shape=circle,inner sep=1mm]{};
\path (0,3) node[draw,shape=circle,inner sep=1mm]{};
\path (0,4) node[draw,shape=circle,inner sep=1mm]{};
\path (0,6) node[draw,shape=circle,inner sep=1mm]{};

\path (0,1) node[fill,shape=circle,inner sep=0.7mm]{};
\path (0,2) node[fill,shape=circle,inner sep=0.7mm]{};
\path (0,5) node[fill,shape=circle,inner sep=0.7mm]{};
\path (0,7) node[fill,shape=circle,inner sep=0.7mm]{};
\path (0,8) node[fill,shape=circle,inner sep=0.7mm]{};
\path (0,9) node[fill,shape=circle,inner sep=0.7mm]{};
\path (0,10) node[fill,shape=circle,inner sep=0.7mm]{};

\path (0.5,0.5) node[draw,shape=circle,inner sep=1mm]{};
\path (0.5,3.5) node[draw,shape=circle,inner sep=1mm]{};
\path (0.5,5.5) node[draw,shape=circle,inner sep=1mm]{};
\path (0.5,8.5) node[draw,shape=circle,inner sep=1mm]{};

\path (0.5,1.5) node[fill,shape=circle,inner sep=0.7mm]{};
\path (0.5,2.5) node[fill,shape=circle,inner sep=0.7mm]{};
\path (0.5,4.5) node[fill,shape=circle,inner sep=0.7mm]{};
\path (0.5,6.5) node[fill,shape=circle,inner sep=0.7mm]{};
\path (0.5,7.5) node[fill,shape=circle,inner sep=0.7mm]{};
\path (0.5,9.5) node[fill,shape=circle,inner sep=0.7mm]{};
\path (0.5,10.5) node[fill,shape=circle,inner sep=0.7mm]{};

\path (1,0) node[draw,shape=circle,inner sep=1mm]{};
\path (1,3) node[draw,shape=circle,inner sep=1mm]{};
\path (1,5) node[draw,shape=circle,inner sep=1mm]{};
\path (1,7) node[draw,shape=circle,inner sep=1mm]{};
\path (1,10) node[draw,shape=circle,inner sep=1mm]{};

\path (1.5,1.5) node[fill,shape=circle,inner sep=0.7mm]{};
\path (1.5,2.5) node[fill,shape=circle,inner sep=0.7mm]{};
\path (1.5,4.5) node[fill,shape=circle,inner sep=0.7mm]{};
\path (1.5,6.5) node[fill,shape=circle,inner sep=0.7mm]{};
\path (1.5,8.5) node[fill,shape=circle,inner sep=0.7mm]{};
\path (1.5,9.5) node[fill,shape=circle,inner sep=0.7mm]{};
\path (1.5,10.5) node[fill,shape=circle,inner sep=0.7mm]{};

\path (1,1) node[fill,shape=circle,inner sep=0.7mm]{};
\path (1,2) node[fill,shape=circle,inner sep=0.7mm]{};
\path (1,4) node[fill,shape=circle,inner sep=0.7mm]{};
\path (1,6) node[fill,shape=circle,inner sep=0.7mm]{};
\path (1,8) node[fill,shape=circle,inner sep=0.7mm]{};
\path (1,9) node[fill,shape=circle,inner sep=0.7mm]{};
\path (1,11) node[fill,shape=circle,inner sep=0.7mm]{};

\path (1.5,0.5) node[draw,shape=circle,inner sep=1mm]{};
\path (1.5,3.5) node[draw,shape=circle,inner sep=1mm]{};
\path (1.5,5.5) node[draw,shape=circle,inner sep=1mm]{};
\path (1.5,7.5) node[draw,shape=circle,inner sep=1mm]{};

\path (2,0) node[draw,shape=circle,inner sep=1mm]{};
\path (2,2) node[draw,shape=circle,inner sep=1mm]{};
\path (2,4) node[draw,shape=circle,inner sep=1mm]{};
\path (2,6) node[draw,shape=circle,inner sep=1mm]{};
\path (2,8) node[draw,shape=circle,inner sep=1mm]{};

\path (2.5,0.5) node[draw,shape=circle,inner sep=1mm]{};
\path (2.5,3.5) node[draw,shape=circle,inner sep=1mm]{};
\path (2.5,4.5) node[draw,shape=circle,inner sep=1mm]{};
\path (2.5,6.5) node[draw,shape=circle,inner sep=1mm]{};

\path (2,1) node[fill,shape=circle,inner sep=0.7mm]{};
\path (2,3) node[fill,shape=circle,inner sep=0.7mm]{};
\path (2,5) node[fill,shape=circle,inner sep=0.7mm]{};
\path (2,7) node[fill,shape=circle,inner sep=0.7mm]{};
\path (2,9) node[fill,shape=circle,inner sep=0.7mm]{};
\path (2,10) node[fill,shape=circle,inner sep=0.7mm]{};
\path (2,11) node[fill,shape=circle,inner sep=0.7mm]{};

\path (2.5,1.5) node[fill,shape=circle,inner sep=0.7mm]{};
\path (2.5,2.5) node[fill,shape=circle,inner sep=0.7mm]{};
\path (2.5,5.5) node[fill,shape=circle,inner sep=0.7mm]{};
\path (2.5,7.5) node[fill,shape=circle,inner sep=0.7mm]{};
\path (2.5,8.5) node[fill,shape=circle,inner sep=0.7mm]{};
\path (2.5,9.5) node[fill,shape=circle,inner sep=0.7mm]{};
\path (2.5,10.5) node[fill,shape=circle,inner sep=0.7mm]{};

\path (-5,5) node[shape=circle,inner sep=0.3mm](valley1){valley: $q^2t$};
\draw[dashed] (A02) -- (valley1);
\draw[dashed] (A06) -- (valley1);
\path (-2,4.5) node[shape=circle,inner sep=0.3mm](peak){peak: $-qt$};
\draw[dashed] (B02) -- (peak);
\draw[dashed] (B05) -- (peak);

\end{scope}


\end{tikzpicture}
\qquad \qquad
\includegraphics[height=3.6in,width=1.4in]{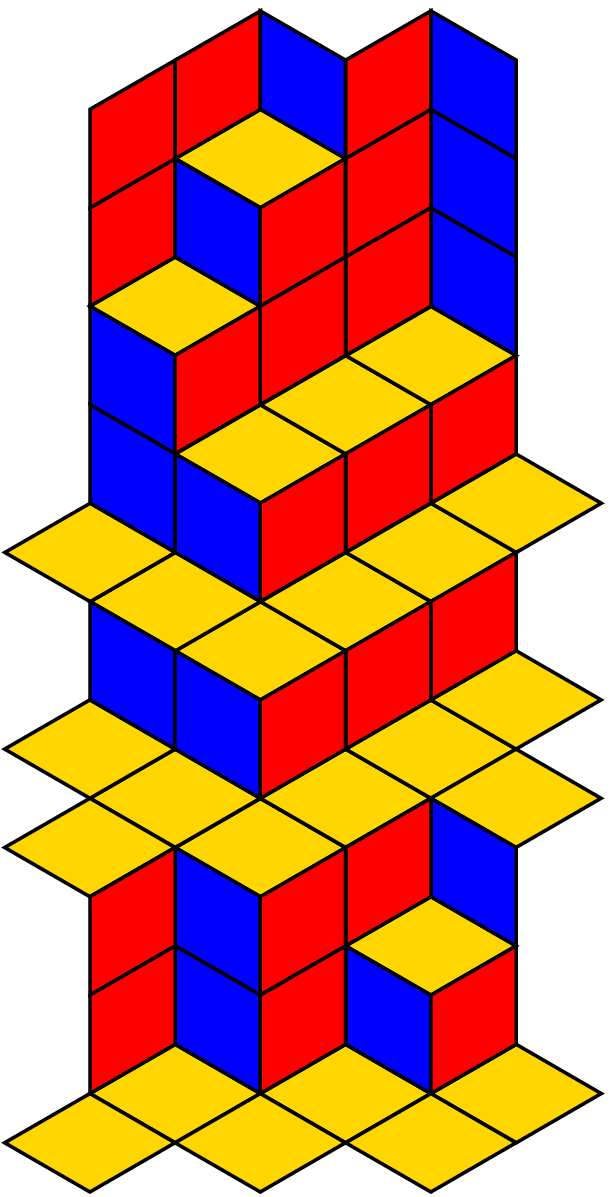}
\end{center}

\begin{equation*} 
 \mathfrak{c} = ((3,2,2), (5,3,2), (6,4,3,2), (4,3,2), (4,3,2,1), (3,2,2)) 
\end{equation*}

The bijection between the path model and the tiling model is clear. The white vertices correspond to the yellow tiles.
Each upstep of a path corresponds to a red tile. Each downstep of a path corresponds to a blue tile.

The following definitions will come in handy when we wish to prove Theorem \ref{theorem2}:

\begin{definition}
We say that a vertex of the lattice is \emph{occupied} or \emph{black} if there is a path passing through that vertex,
otherwise we say that the vertex is \emph{unoccupied} or \emph{white}
\end{definition}

\begin{definition}
A ``cube'' in the non-intersecting lattice path model corresponds to a pair of vertices $u = (x,y_1)$ and $v = (x,y_2)$
with $u$ coloured black, $v$ coloured white and $y_1 < y_2$. 
\end{definition}

\begin{definition}
 A \emph{surface cube} is a cube $(u,v)$ such that if $u = (x,y_1)$ and $v = (x,y_2)$
then for all $w = (x,y')$ with $y_1 < y' < y_2$ the vertex $w$ is coloured white. 
\end{definition}

\begin{definition}
 The \emph{level} of a surface cube $(u,v)$ is $y_2 - y_1$ where $u = (x,y_1)$ and $v = (x,y_2)$.
\end{definition}

Observe that surface cubes are naturally in bijection with the yellow tiles in the rhombus tiling model.

\begin{definition}
The \emph{path associated to the cube $(u,v)$} is the path which passes through the black vertex $v$.
\end{definition}

\begin{definition}
A \emph{valley} cube is a cube $(u,v)$ for which the associated path takes a down step just before passing through $v$,
followed immediately by an upstep.
\end{definition}

\begin{definition}
A \emph{peak} cube is a cube $(u,v)$ for which the associated path takes an up step just before passing through $v$,
followed immediately by a downstep.
\end{definition}

We have marked one peak cube and one value cube on the diagram, together with their contribution to the alphabet $\mathfrak{D}_\mathfrak{c}(q,t)$ in Theorem \ref{comb-weight}. Since we are working on a cylinder, the first vertical is identified with the last vertical in such a way that each path forms a closed loop.

We may mention at this point that in the Hall--Littlewood case we have $q=0$, and the only boxes which contribute to the sum in $\mathfrak{D}_\mathfrak{c}(q,t)$ are those with arm-length zero. Since there is a bijection between such cubes and the ``yellow'' tiles of the rhombus tiling model, and since the leg length of the cube is precisely the \emph{level}, as indicated in Figure 2 of \cite{vuletic2} and Figure 3 of \cite{vuletic}, of the border strip in which the corresponding tile lies, it follows that Theorem \ref{comb-weight} reduces the combinatorial formula for the Hall--Littlewood weight function in the plane partition case in \cite{vuletic2} and both the reverse plane partition and the cylindric plane partition case in \cite{vuletic}.

\hfill \break

After this slight digression, we shall now describe explicitly the bijection between cylindric plane partitions and families of
non-intersecting lattice paths on the cylinder.
The impatient reader may wish to skip the technical details in the section, and simply remark that the parts of the individual
partitions in the interlacing sequence model may be read of the ``heights'' of the corresponding surface cubes
in the rhombus tiling model.

\begin{definition}
The \emph{cylinder of period $T$} is the triangular lattice with vertices $(x,y)$ where either both $x$ and $y$ are even or both $x$ and $y$ are odd, and for which $0 \leq x \leq T$.
\end{definition}

In our example we have $T = 5$. We have drawn a segment of the cylinder corresponding to $0 \leq y \leq 22$.

\begin{definition}
A \emph{path} on the cylinder of period $T$ is a sequence of integers $(y_0, y_1, \ldots, y_T)$ with $y_0$ even such that for each $k$ we have either $y_{k+1} = y_k +1$ or $y_{k+1} = y_k-1$.
\end{definition}


\begin{lemma}\label{binpaths}
Each path on the cylinder of period $T$ may be uniquely
encoded by its starting position $y_0$ and the binary string $p$
given by $p_k = 1$ if $y_{k+1} = y_k + 1$ and $p_k = 0$ otherwise.
\end{lemma}

\begin{definition}
A \emph{family of non-intersecting lattice paths} on the cylinder of period $T$ is a collection of paths:
\begin{align*}
p_1 & = (y_0^1, y_1^1, \ldots y_T^1) \\
p_2 & = (y_0^2, y_1^2, \ldots y_T^2) \\
 & \cdots \\
p_m & = (y_0^m, y_1^m, \ldots y_T^m)
\end{align*}
satisfying $y_{k}^{i+1} > y_k^i$ for all $k$ and $i$ as well as 
$y_0^{i+1} - y_0^i = y_T^{i+1} - y_T^{i}$ for all $i$.
\end{definition}

Note that the second condition is necessary in order to ensure that it is possible to take the cylindric quotient identifying the vertices $(0,y)$ with the vertices $(T,y+d)$ for $d = m-n$ where $m$ is the number of ones in the profile, and $n$ is the number of zeros.

The paths in our example may be encoded using Lemma \ref{binpaths} as follows:
\begin{align*}
p_1 & =(1,0,1,0,1), \qquad y_0^1 = 2 \\
p_2 & =(1,0,1,1,0), \qquad y_0^2 = 4 \\
p_3 & =(0,0,1,1,1), \qquad y_0^3 = 10 \\
p_4 & =(0,0,1,1,1), \qquad y_0^4 = 14 \\
p_5 & =(0,1,1,1,0), \qquad y_0^5 = 16 \\
p_6 & =(1,0,1,1,0), \qquad y_0^6 = 18 \\
p_7 & =(1,1,0,1,0), \qquad y_0^7 = 20 \\
\end{align*}

\begin{definition}
A family of non-intersecting lattice paths on the cylinder of period $T$ is said to be \emph{minimal with $m$ paths} if there is some $i$ such that $y_m^i - y_{m-1}^i > 2$.
\end{definition}

\begin{definition}
The \emph{profile} of a minimal family of $m$ non-intersecting lattice paths is the binary string associated to the $m$th path.
\end{definition}

Our example family of non-intersecting lattice paths is minimal with $7$ paths. Its profile is $\pi = 11010$.
The notion of minimal versus non-minimal cylindric plane partitions is roughly analogous to the notion of minimal versus non-minimal profiles for partitions
(see section \ref{definitions}).

\begin{definition}
The \emph{vertical reading} of a minimal family of $m$ non-intersecting
lattice paths is the sequence of binary strings $\rho_0, \rho_1, \ldots \rho_T$ obtained by reading, for each $k$, vertically upwards from the vertex $(k,y^1_k)$ to the vertex $(k,y^m_k)$, and recording a $0$ each time the vertex is occupied and a $1$ each time the vertex is unoccupied. 
\end{definition}

The vertical reading of our example family of non-intersecting lattice paths is the following:
\begin{align*}
\rho^0 & = 110010111 \quad \quad \rho^3 = 110101011 \\
\rho^1 & = 110101101 \quad \quad \rho^4 = 1010101011 \\
\rho^2 & = 110101011 \quad \quad \rho^5 = 110010111
\end{align*}
Observe that $\rho_0 = \rho_T$.

\begin{proposition}
Let $(\rho_0, \rho_1, \ldots, \rho_T)$ be a sequence of binary strings
arising from the vertical reading of a minimal family of non-intersecting lattice paths on a cylinder of period $T$ and profile $\pi$. For each $k \in \{1,2, \ldots T\}$ let $\mu^k$ denote the partition whose profile is 
given by $\rho_k$. If $\pi_k = 1$ then $\mu^{k} / \mu^{k-1}$ is a horizontal strip, otherwise if $\pi_k = 0$ then $\mu^{k-1} / \mu^{k}$ is a horizontal strip.
\begin{proof}
Follows immediately from Lemma \ref{hstripbin}
\end{proof}
\end{proposition}

Let $\mathfrak{c} = (\mu^0, \mu^1, \ldots, \mu^T)$ be an arbitrary cylindric plane partition with profile $\pi$.
For each $i$ let us define: 
\begin{equation}\label{horizontal-paths}
 \boxed{p_i(\mathfrak{c}) = ((\mu^0)'_i - (\mu^1)'_i + \pi_0, \ldots,
 (\mu^{T})'_i - (\mu^{T-1})'_i + \pi_{T-1}) }
\end{equation}

Note that $p_i(\mathfrak{c})$ encodes information about the length of the $i$-th column of the successive partitions in the interlacing sequence of $\mathfrak{c}$.

\begin{proposition}
For each $i$ we have that $p_i(\mathfrak{c})$ as defined in Equation (\ref{horizontal-paths}) is a binary string.
\begin{proof}
 If $\pi_k = 1$ then from definition \ref{cylindric-def} it follows that $\mu^k / \mu^{k-1}$ is a horizontal strip and thus by Lemma \ref{hstripcol}
we have that $(\mu^{k-1})' - (\mu^k)' \in \{-1,0\}$.
Similarly, if $\pi_k = 0$ then $\mu^{k-1} / \mu^k$ is a horizontal strip and $(\mu^{k-1})' - (\mu^k)' \in \{1,0\}$
\end{proof}
\end{proposition}

\begin{theorem}\label{path-bijection}
For any binary string $\pi$ of length $T$ there is a bijection between minimal families of non-intersecting lattice paths on the cylinder of period $T$ with profile $\pi$, and cylindric plane partitions of profile $\pi$. 
\begin{proof}
The map from families of non-intersecting lattice paths on the cylinder to interlacing sequences is given by taking vertical readings, and then translating from profiles to partitions. Conversely, the family of non-intersecting lattice paths
associated to a given interlacing sequence $\mathfrak{c}$ is given by
\[ \{(p_1(\mathfrak{c}),\tau_1), (p_2(\mathfrak{c}), \tau_2), \ldots (p_m(\mathfrak{c}), \tau_m)\}\]
where $\tau_i$ is the position of $i$-th one in the profile of $\mu^0$ and $p_i(\mathfrak{c})$ is defined in Equation \ref{horizontal-paths}.
\end{proof}
\end{theorem}

\section{Symmetric Functions and Macdonald Polynomials}\label{symfun}
Let $\Lambda_{q,t}$ denote the ring of symmetric functions over the field of rational functions in the indeterminants $q$ and $t$.
Whenever possible we shall suppress in our notation any mention to the variables in which the functions are symmetric.
When we must mention the variables explicitly we shall make use of the \emph{plethystic notation} \cite{garsia}.

In the plethystic notation addition corresponds to the union of two sets and multiplication corresponds to the cartesian product. 
For example, we write:
\begin{equation} X = x_1 + x_2 + \cdots \end{equation}
to denote the set of variables $\{x_1, x_2, \ldots \}$. We also write:
\begin{equation} XY = (x_1 + x_2 + \cdots)(y_1 + y_2, \ldots) \end{equation}
to denote the set of variables $\{x_1 y_1, x_1 y_2, \ldots, x_2 y_1, \ldots x_2 y_2 \ldots \}$.

Let us denote the generating function for the complete symmetric functions by:
\begin{equation} 
 \Omega[Xz] = \prod_i \frac{1}{1-x_iz} = \sum_n h_n z^n 
\end{equation}

The \emph{plethystic negation} of an alphabet may be defined as follows:
\begin{equation}
  \Omega[-Xz] = \prod_i (1-x_iz) = \sum_n (-1)^n e_n z^n 
\end{equation}

where the $e_n$ are the elementary symmetric functions.
In particular, as an alphabet the expression $\frac{1-t}{1-q}$ is to be interpreted as the alphabet:
\begin{equation}
 \{ 1, q, q^2, q^3, \cdots, -t, -tq, -tq^2, -tq^3. \cdots \} 
\end{equation}
So that, for example:
\begin{equation} 
\Omega \left [ \frac{1-t}{1-q}\right ] = \prod_n \frac{1 - tq^{n-1}}{1-q^n} = \frac{(t;q)_\infty}{(q;q)_\infty} 
\end{equation}

We shall define:
\begin{equation} 
 \Omega_{q,t}[Xz] = \Omega \left [ \frac{1-t}{1-q}Xz \right ] = \prod_i \frac{(tx_iz;q)_\infty}{(qx_iz;q)_\infty} 
\end{equation}

Recall that the \emph{Schur functions} are an orthonormal basis for $\Lambda$ with respect to the \emph{Hall Inner product}
\begin{equation} 
 \langle S_\lambda \,|\, S_\mu \rangle = \delta_{\lambda, \mu} 
\end{equation}

The \emph{Cauchy Kernel} is given by:
\begin{equation}
  \Omega[XY] = \prod_{i,j} \frac{1}{1-x_i y_j} = \sum_\lambda S_\lambda(X) S_\lambda(Y) 
\end{equation}

Let $\langle -\,|\,-\rangle_{q,t}$ denote the inner product associated to the $(q,t)$-deformed Cauchy Kernel: 
\begin{equation}
 \Omega_{q,t}[XY] = \Omega \left [ XY \frac{1-t}{1-q} \right ] = \prod_{i,j} \frac{(tx_i y_j;q)_\infty}{(x_i y_j;q)_\infty} 
\end{equation}

where as usual:
\begin{equation} 
(a;q)_n = \prod_{i = 0}^n (1-aq^i) 
\end{equation}

The operator $\Omega^*_{q,t}[Xz]$ is defined to be adjoint to the operator $\Omega_{q,t}[Xz]$ with respect to the Macdonald inner product.

\begin{equation*}
 \langle f(X) \, | \, \Omega^*_{q,t}[Xz] g(X) \rangle_{q,t} = \langle \Omega_{q,t}[Xz] f(X) \, | \, g(X) \rangle_{q,t}  
\end{equation*}

The Macdonald polynomials $\{P_\lambda(X;q,t)\}$ are an orthogonal (but not orthonormal) basis for $\Lambda_{q,t}$ with respect to this inner product
(\cite{macdonald} page 338 -- 340):

\begin{equation}
\langle P_\lambda(X;q,t), P_\mu(X;q,t) \rangle_{q,t} = \delta_{\lambda, \mu} b_\lambda(q,t)
\end{equation}

where:

\begin{equation}
b_\lambda(q,t) = \prod_{s \in \lambda} \frac{(1 - q^{a_\lambda(s)+1} t^{\ell_\lambda(s)})}{(1 - q^{a_\lambda(s)} t^{\ell_\lambda(s)+1}) } 
\end{equation}

Here $a_\lambda(i,j)= \lambda_i - j$ denotes the \emph{arm length} of the box $s = (i,j)$ with respect to the partition $\lambda$ 
and $\ell_\lambda(i,j) = \lambda'_j - i$ denotes the \emph{leg length}.

The dual basis is denoted by $\{Q_\lambda(X;q,t)\}$. 

\begin{equation}
 Q_\lambda(X;q,t) = \prod_{s \in \lambda} \frac{1}{b_\lambda(q,t) } P_\lambda(X;q,t)
\end{equation}

The Pieri formulae for the Macdonald polynomials (\cite{macdonald} (page 340 -- 341) may be expressed in the form:
\begin{equation} \label{pieri}
 \Omega[Xz]_{q,t} \, P_\mu(X;q,t) = \sum_{\lambda \in U(\mu)}\psi_{\lambda / \mu}(q,t) \, P_\lambda(X;q,t) z^{|\lambda|-|\mu|}
\end{equation}
\begin{equation} \label{pieri2}
 \Omega^*[Xz]_{q,t} \, P_\lambda(X;q,t) = \sum_{\mu \in D(\lambda)}\varphi_{\lambda / \mu}(q,t) \, P_\mu(X;
q,t) z^{|\lambda|-|\mu|}
\end{equation}
where $U(\mu)$ is the set of partitions which can be obtained from $\mu$ by adding a horizontal strip, 
$D(\lambda)$ denotes the set of all partitions which can be obtained from $\mu$ by removing a horizontal strip,
and $\varphi_{\lambda / \mu}(q,t)$ and $\psi_{\lambda / \mu}(q,t)$ are given in equations \ref{pieri-coeff} and \ref{pieri-coeff2} respectively.

The following two lemmas are essentially due to Garsia, Haiman and Tesla \cite{garsia}. 
They constitute a $(q,t)$-analog of the commutation relations for ``vertex operators'' to be found in Jimbo and Miwa \cite{soliton}

\begin{lemma}
\begin{equation} \Omega^*_{q,t}[Xz] \, P_\lambda(X;q,t) = P_\lambda(X+z;q,t) \end{equation}
\begin{proof}
Let $\{ Q_\lambda(X;q,t) \}$ denote the dual basis to the $\{P_\lambda(X;q,t)\}$ with respect to the Macdonald inner product.
We have:
\begin{align*}
\Omega^*_{q,t}[Xz] \, P_\lambda(X;q,t) 
& = \langle \Omega^*_{q,t}[Yz] \, P_\lambda(Y;q,t)\,|\, \Omega_{q,t}[XY] \rangle_{q,t} \\
& = \langle  P_\lambda(Y;q,t)\,|\, \Omega_{q,t}[Yz] \, \Omega_{q,t}[XY] \rangle_{q,t} \\
& = \langle  P_\lambda(Y;q,t)\,|\, \Omega_{q,t}[(X+z)Y] \rangle_{q,t} \\
& = P_\lambda(X+z;q,t)
\end{align*}
\end{proof}
\end{lemma}

\begin{lemma}\label{commutation}
\begin{equation} 
\Omega_{q,t}^*[Xu] \, \Omega_{q,t}[Xv] = \frac{(tuv;q)_\infty}{(quv;q)_\infty} \, \Omega_{q,t}[Xv] \, \Omega_{q,t}^*[Xu]
\end{equation}
\begin{proof}
\begin{align*}
\Omega_{q,t}^*[Xu] \, \Omega_{q,t}[Xv] \, P_\lambda(X;q,t) 
& = \Omega_{q,t}[(X+u)z] \, P_\lambda(X + u;q,t)\\
& = \Omega_{q,t}[uz] \, \Omega_{q,t}[Xz] \, \Omega^*_{q,t}[Xz] \, P_\lambda(X;q,t) \\
& = \prod_{n \geq 0} \frac{(tuv;q)_\infty}{(quv;q)_\infty} \, \Omega_{q,t} [Xz] \, \Omega^*_{q,t}[Xz] \, P_\lambda(X;q,t)
\end{align*}
\end{proof}
\end{lemma}

\section{Proof of theorem \ref{qt-borodin}}\label{theorem1}

In this section we prove a $(q,t)$-analog of Borodin's formula for the enumeration of cylindric plane partitions (Theorem \ref{qt-borodin}).

We begin with a number of small lemmas. 
Let $D_z$ denote the `degree'' operator:
\begin{equation} 
D_z P_\lambda[X] = z^{|\lambda|} P_\lambda[X] 
\end{equation}

The degree operator satisfies the following commutation relations:
\begin{lemma} \label{degree}
\begin{align}
D_z \, \Omega_{q,t}[Xu] & = \Omega_{q,t}[Xuz] \, D_z \\
D_z \, \Omega^*_{q,t}[Xu] & = \Omega^*_{q,t}[Xuz^{-1}] \, D_z 
\end{align}
\begin{proof}
This fact follows immediately from the action of $\Omega_{q,t}[Xu]$ and $\Omega_{q,t}[Xv]$ on Macdonald polynomials 
(equations \ref{pieri} and \ref{pieri2}).
\end{proof}
\end{lemma}

For notational convenience we shall define:
\begin{align}
G^0(z) & = \Omega_{q,t}[Xz] \\
G^1(z) & = \Omega^*_{q,t}[Xz]
\end{align}

\begin{lemma}  \label{lhs}
The left hand side of the refined version of equation \ref{maintheorem} may be expressed in the form:
\begin{equation} 
\lhs(\pi) = \sum_\mu \langle Q_\mu \,|\, G^{\pi_0}(u_0) G^{\pi_1}(u_1)  \cdots G^{\pi_T}(u_T) D_w \, P_\mu  \rangle_{q,t} 
\end{equation}
where:
\begin{align} \label{specialization}
w & = z_0 z_1 \cdots z_{T-1} \\ 
u_k & = \begin{cases}
z_0 z_1 \cdots z_{k-1} & \text{ if $\pi_k = 1$} \\
z_0^{-1} z_1^{-1} \cdots z_{k-1}^{-1} & \text{ if $\pi_k = 0$}
\end{cases}\label{specialization2}
\end{align}

\begin{proof}
From the ``interlacing sequence'' definition of a cylindric plane partition it is clear that a cylindric plane partition is constructed
by successively adding and removing horizontal strips. 

The presence of the $(q,t)$-Pieri coefficients in the definition of the weight function (equation \ref{weight})
come directly from the action of the operators $\Omega_{q,t}$ and $\Omega^*_{q,t}$ given in equations \ref{pieri} and \ref{pieri2}.
The degree operator $D_z$ is used to keep track of the number of cubes in the resulting cylindric plane partition.

Using the fact that the Macdonald $P$-functions are orthogonal with respect to the Macdonald $Q$-functions 
we may write:
\begin{equation}
\lhs(\pi) = \sum_\mu \langle Q_\mu \,|\, D_{z_0} \, G^{\pi_0}(1) \, D_{z_1} \, G^{\pi_1}(1)   \cdots D_{z_{T-1}} \, G^{\pi_T}(1) \, P_\mu  \rangle_{q,t} 
\end{equation}

It remains to commute all the shift operators to the right hand side using Lemma \ref{degree}.
\end{proof}
\end{lemma}

Next let us define:
\begin{definition}\label{M-def}
\begin{equation} 
M_\pi(m)  = \sum_\mu \langle Q_\mu \,|\, \prod_{\substack{k =1\\ \pi_k = 0}}^T \Omega_{q,t} [Xu_kw^m] 
\prod_{\substack{k =1 \\ \pi_k = 1}}^T \Omega_{q,t}^*[Xu_k] \, D_w \, P_\mu \rangle_{q,t} 
\end{equation}
\end{definition}

\begin{lemma} \label{cylindric-shift}
\[ M_\pi(m) =  \prod_{\substack{(i,j) \\ \pi_i \neq \pi_j}} \frac{(tu_i u_j w^{m+1};q)_\infty}{(u_i u_j w^{m+1};q)_\infty} \, M_\pi(m+1) \]
\begin{proof}
This is a straightforward calculation. Using the fact that the $\{P_\lambda\}$ are orthogonal to the $\{Q_\lambda\}$ we may write:
\begin{align}
M_\pi(m) & = \sum_{\mu, \lambda} 
\langle Q_\mu \,|\, \prod_{\substack{k =1\\ \pi_k = 0}}^T \Omega_{q,t}[Xu_kw^m] \, P_\lambda \rangle_{q,t} 
\langle Q_\lambda \,|\,\prod_{\substack{k =1 \\ \pi_k = 1}}^T \Omega_{q,t}^*[Xu_k] \, D_w \, P_\mu \rangle_{q,t} \\
& = \sum_{\mu, \lambda} 
\langle Q_\lambda \,|\,\prod_{\substack{k =1 \\ \pi_k = 1}}^T \Omega_{q,t}^*[Xu_k] \, D_w \, P_\mu \rangle_{q,t} 
\langle Q_\mu \,|\, \prod_{\substack{k =1\\ \pi_k = 0}}^T \Omega_{q,t}[Xu_kw^m] \, P_\lambda \rangle_{q,t} \\
& = \sum_{\lambda} 
\langle Q_\lambda \,|\,\prod_{\substack{k =1 \\ \pi_k = 1}}^T \Omega_{q,t}^*[Xu_k] \, 
D_w \,\prod_{\substack{k =1\\ \pi_k = 0}}^T \Omega_{q,t}[Xu_kw^m] \, P_\lambda \rangle_{q,t} 
\end{align}
Next applying the commutation relations of Lemma \ref{degree} and Lemma \ref{commutation} we have:
\begin{align}
M_\pi(m) & = \sum_{\lambda} 
\langle Q_\lambda \,|\,\prod_{\substack{k =1 \\ \pi_k = 1}}^T \Omega_{q,t}^*[Xu_k] \, 
D_w \,\prod_{\substack{k =1\\ \pi_k = 0}}^T \Omega_{q,t}[Xu_kw^m] \, P_\lambda \rangle_{q,t}  \\ 
& = \sum_{\lambda} 
\langle Q_\lambda \,|\,\prod_{\substack{k =1 \\ \pi_k = 1}}^T \Omega_{q,t}^*[Xu_k] \, 
\,\prod_{\substack{k =1\\ \pi_k = 0}}^T \Omega_{q,t}[Xu_kw^{m+1}] \, D_w \, P_\lambda \rangle_{q,t} \\
& = \prod_{\substack{(i,j) \\ \pi_i \neq \pi_j}} \frac{(tu_i u_j w^{m+1};q)_\infty}{(u_i u_j w^{m+1};q)_\infty} \sum_{\lambda} 
\langle Q_\lambda \,| \,\prod_{\substack{k =1\\ \pi_k = 0}}^T \Omega_{q,t}[Xu_kw^{m+1}] \, 
\prod_{\substack{k =1 \\ \pi_k = 1}}^T \Omega_{q,t}^*[Xu_k] \, D_w \, P_\lambda \rangle_{q,t} \\
& =  \prod_{\substack{(i,j) \\ \pi_i \neq \pi_j}} \frac{(tu_i u_j w^{m+1};q)_\infty}{(u_i u_j w^{m+1};q)_\infty} \, M_\pi(m+1)
\end{align}
\end{proof}
\end{lemma}
In the limit we have:

\begin{lemma} \label{infinity}
\begin{equation} 
M_\pi(\infty)  = \prod_{n \geq 1} \frac{1}{1-w^n} 
\end{equation}
\begin{proof}
In order for this limit to even make sense, we must have $|z_i| < 1$ for all $i$, in which case:
\[ \lim_{m \to \infty} \Omega_{q,t}[X u_k \omega^m] = 1 \]

Since $\Omega_{q,t}^*[X u_k]$ is a degree lowering operator, it follows that:
\begin{align*}
 \lim_{m \to \infty} M_\pi(m) 
& =  \sum_\mu \langle Q_\mu \,|\,  \prod_{\substack{k =1 \\ \pi_k = 1}}^T \Omega_{q,t}^*[Xu_k] \, D_w \, P_\mu \rangle_{q,t} \\
& = \sum_\mu \langle Q_\mu | D_w \, P_\mu \rangle_{q,t} \\
& = \sum_\mu \omega^{|\mu|} \\
& = \prod_{n \geq 1} \frac{1}{1-w^n} 
\end{align*}

\end{proof}
\end{lemma}


The proof of the refined version of Theorem \ref{qt-borodin} now proceeds as follows. We begin by applying Lemma \ref{lhs}

\begin{align*}
\phantom{=} \sum_{\mathfrak{c} \in \CPP(\pi)} W_\mathfrak{c}(q,t) z^{|\mathfrak{p}|} 
& = \sum_\mu \langle Q_\mu \,|\, G^{\pi_0}(u_0) G^{\pi_1}(u_1)  \cdots G^{\pi_T}(u_T) D_w \, P_\mu  \rangle_{q,t}  \\
\end{align*}
Next we repeatedly applies the commutation relations of Lemma \ref{commutation}, followed by definition \ref{M-def}.
\begin{align*}
& = \prod_{\substack{i < j \\ \pi_i > \pi_j} }  \frac{(t u_i u_j;q)_\infty}{(u_i u_j;q)_\infty}
\sum_\mu \langle Q_\mu \,|\, \prod_{\substack{k =1\\ \pi_k = 0}}^T \Omega_{q,t} [Xu_k] 
\prod_{\substack{k =1 \\ \pi_k = 1}}^T \Omega_{q,t}^*[Xu_k] \, D_w \, P_\mu \rangle_{q,t} \\
& = \prod_{\substack{i < j \\ \pi_i > \pi_j} } \frac{(t u_i u_j;q)_\infty}{(u_i u_j;q)_\infty} M_\pi(0) 
\end{align*}
We then repeatedly apply Lemma \ref{cylindric-shift}.
\begin{align*}
& = \prod_{\substack{i < j \\ \pi_i > \pi_j} } \frac{(t u_i u_j;q)_\infty}{(u_i u_j;q)_\infty} \prod_{m \geq 0} 
\left ( \prod_{\substack{(i,j) \\ \pi_i \neq \pi_j}} \frac{(tu_i u_j w^{m+1};q)_\infty}{(u_i u_j w^{m+1};q)_\infty} \right ) M_\pi (\infty) 
\end{align*}
Splitting the second product into two, and combining it with the first we have:

\begin{align*}
& = \prod_{m \geq 1}
\left ( \prod_{\substack{i < j \\ \pi_i > \pi_j} } \frac{(t u_i u_j w^{m-1};q)_\infty}{(u_i u_jw^{m-1};q)_\infty} \right )
\left ( \prod_{\substack{i > j \\ \pi_i > \pi_j} } \frac{(t u_i u_j w^m;q)_\infty}{(u_i u_jw^m;q)_\infty}\right ) M_\pi (\infty)
\end{align*}
Finally, applying Lemma \ref{infinity} we have:
\begin{align*}
& = \prod_{m \geq 1} \frac{1}{1-w^m}
\left ( \prod_{\substack{i < j \\ \pi_i > \pi_j} } \frac{(t u_i u_j w^{m-1};q)_\infty}{(u_i u_jw^{m-1};q)_\infty} \right )
\left ( \prod_{\substack{i > j \\ \pi_i > \pi_j} } \frac{(t u_i u_j w^m;q)_\infty}{(u_i u_jw^m;q)_\infty}\right ) 
\end{align*}

To obtained the non-refined version of the Theorem, it suffices to take the following specialization of variables on both sides:
\begin{align} \label{specialization}
w & = z^{|T|} \\ 
u_k & = \begin{cases}
z^k & \text{ if $\pi_k = 1$} \\
z^{-k} & \text{ if $\pi_k = 0$}
\end{cases}\label{specialization2}
\end{align}

\section{Proof of Theorem \ref{comb-weight}}\label{theorem2}

We begin by making use of the \emph{plethystic notation} to rewrite the $(q,t)$-Pieri coefficients (equations \ref{pieri-coeff} and \ref{pieri-coeff2}) in the following form:
\begin{align}
\varphi_{\lambda / \mu}(q,t) & = \Omega \left [ (q-t) (\mathcal{A}_{\lambda / \mu}(q,t) - \mathcal{B}_{\lambda / \mu}(q,t)) \right ] \\
\psi_{\lambda / \mu}(q,t) & = \Omega \left [ (q-t) (\mathcal{B}'_{\lambda / \mu}(q,t) - \mathcal{A}'_{\lambda / \mu}(q,t)) \right ]
\end{align}
where:
\begin{align}
\mathcal{A}_{\lambda / \mu}(q,t) & = \sum_{s \in C_{\lambda / \mu}} q^{a_\lambda(s)} t^{\ell_\lambda(s)} \\
\mathcal{B}_{\lambda / \mu}(q,t) & = \sum_{s \in \overline{C}_{\lambda / \mu}} q^{a_\mu(s)} t^{\ell_\mu(s)} \\
\mathcal{A}'_{\lambda / \mu}(q,t) & = \sum_{s \not\in C_{\lambda / \mu}} q^{a_\lambda(s)} t^{\ell_\lambda(s)} \\
\mathcal{B}'_{\lambda / \mu}(q,t) & = \sum_{s \not\in \overline{C}_{\lambda / \mu}} q^{a_\mu(s)} t^{\ell_\mu(s)} 
\end{align}

Here we have changed our notation slightly from that used in Macdonald, so that now $C_{\lambda / \mu}$ denotes the set of boxes $s = (i,j) \in \lambda$ such that $\lambda'_j > \mu'_j$ while $\overline{C}_{\lambda / \mu}$ denotes the set of boxes $s = (i,j) \in \mu$ such that $\lambda'_j > \mu'_j$.

Making use of this notation we may rewrite equation \ref{weight} as:
\begin{equation} 
W_{\mathfrak{c}}(q,t) = \Omega \bigl [ (q-t) \mathcal{D}_{\mathfrak{c}}(q,t) \bigr ] 
\end{equation}
where:
\begin{equation}\label{unseparated}
\mathcal{D}_{\mathfrak{c}}(q,t) = \sum_{\substack{k =1 \\ \pi_k = 1}}^T  (\mathcal{A}_{k / {k-1}}- \mathcal{B}_{k / {k-1}}) 
+ \sum_{\substack{k = 1\\ \pi_k = 0}}^T (\mathcal{B}'_{k-1 / {k}} - \mathcal{A}'_{k-1 / {k}}) \\
\end{equation}
To avoid unnecessary indices, we use the convention that: 
\begin{equation} 
\mathcal{X}_{k / {k-1}} = \mathcal{X}_{\mu^k / \mu^{k-1}}(q,t)
\end{equation}

Our goal is to find a simplified expression for $\mathcal{D}_{\mathfrak{c}}(q,t)$. Recall from section \ref{definitions} that in the ``interlacing sequence'' model, a \emph{cube} of the cylindric plane partition $\mathfrak{c}$ corresponds to a \emph{box} of one of the underlying partitions $\mu^k$.

Observe now that each \emph{box} $s \in \mu^k$ contributes to at most two terms in equation \ref{unseparated}, one involving the pair of partitions $\mu^k$ and $\mu^{k-1}$, the other involving the pair of partitions $\mu^k$ and $\mu^{k+1}$.

Regrouping terms, and setting $\pi_{T+1} = \pi_1$ as well as $\mu^{T+1} = \mu^1$ we may write:
\begin{equation}
\mathcal{D}_{\mathfrak{c}}(q,t) = 
 \sum_{\substack{k =1 \\ \pi_k = 1 \\\pi_{k + 1} = 1}}^{T} \mathcal{E}^k_{11}(\mathfrak{c}) 
+ \sum_{\substack{k =1 \\ \pi_k = 0 \\ \pi_{k+1} = 0}}^{T} \mathcal{E}^k_{00}(\mathfrak{c}) 
+ \sum_{\substack{k =1 \\ \pi_k = 0 \\ \pi_{k+1}=1} }^{T} \mathcal{E}^k_{01}(\mathfrak{c}) 
+ \sum_{\substack{k =1 \\ \pi_k =1 \\ \pi_{k+1} = 0}}^{T} \mathcal{E}^k_{10}(\mathfrak{c}) 
\end{equation}

where:
\begin{align}
\mathcal{E}^k_{11}(\mathfrak{c}) & = \mathcal{A}_{k/k-1} - \mathcal{B}_{k+1/k}\\
\mathcal{E}^k_{00}(\mathfrak{c}) & = \mathcal{B}'_{k-1 / k} - \mathcal{A}'_{k/k+1}\\
\mathcal{E}^k_{01}(\mathfrak{c}) & = \mathcal{B}'_{k-1/k} - \mathcal{B}_{k+1/k} \\
\mathcal{E}^k_{10}(\mathfrak{c}) & = \mathcal{A}_{k/k-1}- \mathcal{A}'_{k/k+1}
\end{align}

For each $k$, there is only one term of the form $\mathcal{E}^k_{rs}(\mathfrak{c})$ appearing in the expression for $\mathcal{D}_{\mathfrak{c}}(q,t)$, and this term groups together all contributions from the boxes $s \in \mu^k$.

The next step is to observe that we have a large number of cancellations. For example:
\begin{align*}
\mathcal{E}^k_{11}(\mathfrak{c}) & = (\mathcal{A}_{k / {k-1}} - \mathcal{B}_{{k+1} / {k}})\\
& = \sum_{s \in C_{k / k-1}} q^{a_k(s)} t^{\ell_k(s)} - \sum_{s \in \overline{C}_{k+1/k}} q^{a_k(s)} t^{\ell_k(s)} \\
& = \sum_{s \in \mu^k} \sign_{11}(s) \, q^{a_k(s)} t^{\ell_k(s)}
\end{align*}
where:
\begin{equation} 
\sign_{11}(s) = 
\begin{cases}
1 & \text{if } s \in C_{k / k-1} \text{ and } s \not\in \overline{C}_{k+1/k} \\
0 & \text{if } s \in C_{k / k-1} \text{ and } s \in \overline{C}_{k+1/k} \\
0 & \text{if } s \not\in C_{k / k-1} \text{ and } s \not\in \overline{C}_{k+1/k}\\
-1 & \text{if } s \not\in C_{k / k-1} \text{ and } s \in \overline{C}_{k+1/k}
\end{cases}
\end{equation}
Again we are using a simplified notation:
\begin{align*}
a_k(s) & = a_{\mu^k}(s) \\ 
\ell_k(s) & = \ell_{\mu^k}(s)
\end{align*}
Similarly:
\begin{align*}
\mathcal{E}^k_{00} & = (\mathcal{B}'_{k-1 / k}- \mathcal{A}'_{k/k+1})\\
& = \sum_{s \not\in \overline{C}_{k-1 /k}} q^{a_k(s)} t^{\ell_k(s)} - \sum_{s \not\in C_{k/k+1}} q^{a_k(s)} t^{\ell_k(s)} \\
& = \sum_{s \in \mu^k} \sign_{00}(s) \, q^{a_k(s)} t^{\ell_k(s)}
\end{align*}
where:
\begin{equation} \sign_{00}(s) = 
\begin{cases}
-1 & \text{if } s \in \overline{C}_{k-1 /k} \text{ and } s \not\in C_{k/k+1}\\
0 & \text{if } s \in \overline{C}_{k-1 /k}\text{ and } s \in C_{k/k+1}\\
0 & \text{if } s \not\in \overline{C}_{k-1 /k} \text{ and } s \not\in C_{k/k+1}\\
1 & \text{if } s \not\in \overline{C}_{k-1 /k} \text{ and } s \in C_{k/k+1}
\end{cases}
\end{equation}
Next:
\begin{align*}
\mathcal{E}^k_{01} & = (\mathcal{B}'_{k-1/k} - \mathcal{B}_{k+1/k})\\
& = \sum_{s \not\in \overline{C}_{k-1 / k}} q^{a_k(s)} t^{\ell_k(s)} - \sum_{s \in \overline{C}_{k+1/k}} q^{a_k(s)} t^{\ell_k(s)} \\
& = \sum_{s \in \mu^k} \sign_{01}(s) \, q^{a_k(s)} t^{\ell_k(s)}
\end{align*}
where:
\begin{equation} \sign_{01}(s) = 
\begin{cases}
0 & \text{if } s \in \overline{C}_{k-1 / k} \text{ and } s \not\in \overline{C}_{k+1/k} \\
-1 & \text{if } s \in \overline{C}_{k-1 / k} \text{ and } s \in \overline{C}_{k+1/k} \\
1 & \text{if } s \not\in \overline{C}_{k-1 / k} \text{ and } s \not\in \overline{C}_{k+1/k} \\
0 & \text{if } s \not\in \overline{C}_{k-1 / k} \text{ and } s \in \overline{C}_{k+1/k}(
\end{cases}
\end{equation}
And:
\begin{align*}
\mathcal{E}^k_{10}(\mathfrak{c}) & = (\mathcal{A}_{k/k-1} - \mathcal{A}'_{k/k+1})\\
& = \sum_{s \in C_{k / k-1}} q^{a_k(s)} t^{\ell_k(s)} - \sum_{s \not\in C_{k/k+1}} q^{a_k(s)} t^{\ell_k(s)} \\
& = \sum_{s \in \mu^k} \sign_{10}(s) \, q^{a_k(s)} t^{\ell_k(s)}
\end{align*}
where:
\begin{equation} \sign_{10}(s) = 
\begin{cases}
0 & \text{if } s \in C_{k/ k-1}\text{ and } s \not\in C_{k/k+1} \\
1 & \text{if } s \in C_{k/ k-1} \text{ and } s \in C_{k/k+1} \\
-1 & \text{if } s \not\in C_{k/ k-1} \text{ and } s \not\in C_{k/k+1} \\
0 & \text{if } s \not\in C_{k/ k-1} \text{ and } s \in C_{k/k+1}
\end{cases}
\end{equation}

The final step in the proof is to switch from the ``interlacing sequence'' model of cylindric plane partitions (section \ref{definitions}) to the non-intersecting path model (section \ref{lattice-paths}). This entails a grouping together of all the \emph{cubes} of the cylindric plane partition which belong to the same column of possibly different partitions in the sequence.

Recall that a \emph{cube} in the non-intersecting path model corresponds to a pair of vertices $v_1 = (x,y_1)$ and $v_2 = (x,y_2)$ with $y_1 < y_2$ where $v_1$ is coloured black and $v_2$ is coloured white.

Recall also that the $i$-th path in the non-intersecting path model encodes the length of the $i$-th column in each succeeding partition of the interlacing sequence (see equation (\ref{horizontal-paths}) in section \ref{lattice-paths}).

We shall say that the cube $c = (v_1,v_2)$ is of type $i$ if the black vertex $y_1$ lies on the $i$th path. This is equivalent to saying that the cube $c$ lies in the $i$th column of $\mu^k$ for some $k$.

Now, if $\pi_k = 0$ and $\mu^k = \mu^{k-1}$ then at the $k$th step, all the paths move downwards.
More generally if $\pi_k = 0$ then $\mu^k \preceq \mu^{k-1}$ and the $i$th path moves upwards if and only if the $i$th column of $\mu^k$ is shorter than the corresponding column of $\mu^{k-1}$.

That is to say, at the $k$th step, the $i$th path moves upwards if and only if $c \in \overline{C}_{k-1/k}$ for all $c \in \mu^k$ of type $i$.

If $\pi_k = 1$ and $\mu^k = \mu^{k-1}$ then at the $k$th step, all the paths move upwards.
More generally if $\pi_k = 1$ then $\mu^{k-1} \preceq \mu^{k}$ and the $i$th path moves downwards if and only if the $i$th column of $\mu^k$ is longer than the corresponding column of $\mu^{k-1}$.

That is to say, at the $k$th step, the $i$th path moves downwards if and only if $c \in C_{k/k-1}$ for all $c \in \mu^k$ of type $i$.

In a similar spirit, $\pi_{k+1} = 0$ and $\mu^{k+1} = \mu^{k}$ then at the $(k+1)$th step, all the paths move downwards.
More generally if $\pi_{k+1} = 0$ then $\mu^{k+1} \preceq \mu^{k}$ and the $i$th path moves up if and only if the $i$th column of $\mu^{k+1}$ is shorter than the corresponding column of $\mu^{k}$.

That is to say, at the $(k+1)$th step, the $i$th path moves upwards if and only if $c \in C_{k/k+1}$ for all $c \in \mu^k$ of type $i$.

If $\pi_{k+1} = 0$ and $\mu^{k+1} = \mu^{k}$ then at the $(k+1)$th step, all the paths move upwards.
More generally if $\pi_k = 1$ then $\mu^{k} \preceq \mu^{k+1}$ and the $i$th path moves downwards if and only if the $i$th column of $\mu^{k+1}$ is longer than the corresponding column of $\mu^{k}$.

That is to say, at the $(k+1)$th step, the $i$th path moves downwards if and only if $c \in \overline{C}_{k+1/k}$ for all $c \in \mu^k$ of type $i$.

Last but not least, one may check that the signs agree in all 16 possible cases.

\section{Conclusion} \label{conclusion}
We have proven a $(q,t)$-analog for the enumeration of cylindrical plane partitions which was first discovered 
by Borodin \cite{borodin}. 
The proof relies on certain commutation relations between operators acting on Macdonald polynomials 
which are essentially due to Garsia, Haiman and Tesla \cite{garsia}. In the Schur case, these operators are precisely the \emph{vertex operators} of Jimbo and Miwa \cite{soliton}.

By interpreting cylindric plane partitions as non-intersecting lattice paths we have also given an explicit combinatorial interpretation of the weight function in this $(q,t)$-analog. This greatly simplifies the expression for the weight given by Okada \cite{okada} in the reverse plane partition case. It also reduces to the expression given by Corteel, Savelief and Vuleti{\'c} in the Hall--Littlewood case \cite{vuletic}.

Although we have worked exclusively with symmetric functions, it ought to be possible to reformulate everything in terms of operators acting on fermionic fock space. 
There exists a Hall--Littlewood version of the Boson--Fermion correspondence which is given in Jing \cite{jing}, and which has been used by Foda and Wheeler \cite{wheeler2} to give a fermionic perspective on the $t$-deformed enumeration of plane partitions. 
No such correspondence is known in the Macdonald case. In particular, there is no simple formula to be found in the literature for the multiplication of a Macdonald polynomial by a power sum (Murnaghan--Nakayama lemma). 

The proof of the $(q,t)$-Pieri formula given in Macdonald (\cite{macdonald} pages 331-341) leaves much to be desired in terms of clarity and simplicity, and no simple formula is given for the $(q,t)$-analog of the Littlewood--Richardson rule - as a weighted sum over puzzles, for example \cite{puzzles}.

As a final remark, we mention that the coefficients $\varphi_T(q,t)$ and $\psi_T(q,t)$ given in equations (7.11) and (7.11') of Macdonald (page 346) 
for the expansion of Macdonald polynomials in terms of monomials
may be given an almost identical treatment to our weight function for the enumeration of cylindric plane partitions. 
Even in the Schur case, from the combinatorial definition as a sum over semi-standard Young tableaux it is not completely obvious that the resulting functions are symmetric,
with the additional $(q,t)$-weight, the symmetry is even less obvious in the Macdonald case.

We shall work on these problems and present them in a forthcoming paper.

\bibliographystyle{alpha}
\bibliography{references}

\end{document}